# The Mean Curvature Flow Smoothes Lipschitz Submanifolds

Mu-Tao Wang

February 15, 2002, revised March 20, 2003

**Abstract**

The mean curvature flow is the gradient flow of volume functionals on the space of submanifolds. We prove a fundamental regularity result of the mean curvature flow in this paper: a Lipschitz submanifold with small local Lipschitz norm becomes smooth instantly along the mean curvature flow. This generalizes the regularity theorem of Ecker and Huisken for Lipschitz hypersurfaces. In particular, any submanifold of the Euclidean space with a continuous induced metric can be smoothed out by the mean curvature flow. The smallness assumption is necessary in the higher codimension case in view of an example of Lawson and Osserman. The stationary phase of the mean curvature flow corresponds to minimal submanifolds. Our result thus generalizes Morrey's classical theorem on the smoothness of $C^1$ minimal submanifolds.

## 1  Introduction

Let $F : \Sigma \to \mathbb{R}^N$ be an isometric immersion of a compact $n$-dimensional manifold $\Sigma$ in the Euclidean space. The mean curvature flow of $F$ is a family of immersions $F_t : \Sigma \to \mathbb{R}^N$ that satisfies

$$\frac{d}{dt} F_t(x) = H(x, t) \qquad (1.1)$$
$$F_0 = F$$

where $H(x,t)$ is the mean curvature vector of $\Sigma_t \equiv F_t(\Sigma)$ at $F_t(x)$. In terms of local coordinates $x^1, \cdots, x^n$ on $\Sigma$, the mean curvature flow is the solution

$$F_t = F^A(x^1, \cdots, x^n, t), \ A = 1, \cdots, N$$

of the following system of parabolic equations

$$\frac{\partial F^A}{\partial t} = \sum_{i,j,B} g^{ij} P_B^A \frac{\partial^2 F^B}{\partial x^i \partial x^j}, \ A = 1, \cdots, N$$

where $g^{ij} = (g_{ij})^{-1}$ is the inverse of the induced metric $g_{ij} = \sum_A \frac{\partial F^A}{\partial x^i} \frac{\partial F^A}{\partial x^j}$. $P_B^A = \delta_B^A - g^{kl} \frac{\partial F^A}{\partial x^k} \frac{\partial F^B}{\partial x^l}$ at any point $p \in \Sigma_t$ is the projection to the normal space, the orthogonal complement to the tangent space of $\Sigma_t$ at $p$.

By the first variation formula, the mean curvature flow is the gradient flow of volume functionals on the space of submanifolds. The well-known formula $H = \Delta_\Sigma F$ suggests (1.1) should be viewed as the heat equation of submanifolds. Here $\Delta_\Sigma$ denotes the Laplace operator of the induced metric $g_{ij}$ on $\Sigma$.

This paper discusses the smoothing property of the mean curvature flow. It was proved by Ecker and Huisken in [4] that any complete hypersurface in $\mathbb{R}^N$ satisfying local uniform Lipschitz condition becomes smooth instantly along the mean curvature flow. In this paper, we proved the following generalization in the arbitrary codimension case.

**Theorem A** *Let $\Sigma$ be a compact n-dimensional Lipschitz submanifold of $\mathbb{R}^{n+m}$. There exists a positive constant $K$ depending on $n$ and $m$ such that if $\Sigma$ satisfies the $K$ local Lipschitz condition, then the mean curvature flow of $\Sigma$ has a smooth solution on some time interval $(0, T]$.*

The time $T$ can be estimated in terms of $K$. The $K$ local Lipschitz condition will be defined in §5. Under this condition $\Sigma$ is allowed to have corners but the corners of $\Sigma$ cannot be too sharp. This condition should be necessary in view of an example of Lawson and Osserman. In [5], they constructed a stable minimal cone in $\mathbb{R}^7$ that is a Lipschitz graph over $\mathbb{R}^4$.

Since minimal submanifolds are stationary phase of the mean curvature flow, Theorem A indeed generalizes a classical theorem of Morrey [6] which asserts any $C^1$ minimal submanifold is smooth.

**Corollary A** *Any minimal submanifold satisfying the $K$ local Lipschitz condition is smooth.*

We remark that any $C^1$ submanifold satisfies the $K$ local Lipschitz condition.

The strategy of the proof follows that of [4]. We prove an a priori curvature estimates for smooth submanifolds in terms of a controlled Lipschitz bound. We then approximate $\Sigma$ by smooth submanifolds and show the mean curvature flow of the approximating submanifolds converges to that of $\Sigma$ which is now smooth for $t > 0$. The global version of such estimate in arbitrary codimension already appeared in [7] and [8]. I would like to thank Professor R. Schoen who suggested me to localize these estimates.

I am grateful to Professor D. H. Phong and Professor S.-T. Yau for their constant advice, encouragement and support. I would also like to thank Professor G. Perelman for interesting remarks and Professor T. Ilmanen for helpful discussions.

## 2   Preliminaries

Let $F_t : \Sigma \to \mathbb{R}^{n+m}$ be the mean curvature flow of a compact $n$-dimensional smooth submanifold. We shall assume $\Sigma_t$ is locally given as a graph over some $\mathbb{R}^n$. We choose an orthonormal basis $\{a_i\}_{i=1\cdots n}$ for this $\mathbb{R}^n$ and an orthonormal basis $\{a_\alpha\}_{\alpha=n+1\cdots n+m}$ for the orthogonal complement $\mathbb{R}^m$. Then $\Omega = a_1^* \wedge \cdots \wedge a_n^*$ can be viewed as an $n$ form on $\mathbb{R}^{n+m}$ which corresponds to the volume form of the $\mathbb{R}^n$. Define the function $*\Omega$ on $\Sigma_t$ by

$$*\Omega = \frac{\Omega(\frac{\partial F_t}{\partial x^1}, \cdots \frac{\partial F_t}{\partial x^n})}{\sqrt{\det g_{ij}}}.$$

This function is the Jacobian of the projection from $\Sigma_t$ onto $\mathbb{R}^n$ and $0 < *\Omega \leq 1$ whenever $\Sigma_t$ can be written as a graph over $\mathbb{R}^n$. At any point $p \in \Sigma_t$, we can choose an oriented orthonormal basis $\{e_i\}_{i=1\cdots n}$ for the tangent space $T_p\Sigma_t$ and $*\Omega = \Omega(e_1, \cdots, e_n)$. We can also choose an orthonormal basis $\{e_\alpha\}_{\alpha=n+1\cdots n+m}$ for the normal space $N_p\Sigma_t$, then the second fundamental form $A = \{h_{\alpha ij}\}$ is represented by $h_{\alpha ij} = \langle \nabla_{e_i} e_j, e_\alpha \rangle$. The convention that $i, j, k$ denote tangent indices and $\alpha, \beta, \gamma$ normal indices are followed.

With these notations we recall the following evolution equations from [7] and [8].

**Lemma 2.1** *Let $\Omega$ be a parallel $n$ form on $R^{n+m}$. Along the mean curvature*

*flow, we have*

$$(\frac{d}{dt} - \Delta) * \Omega = *\Omega |A|^2 - 2 \sum_{\alpha,\beta,k} [\Omega_{\alpha\beta 3 \cdots n} h_{\alpha 1 k} h_{\beta 2 k} + \cdots + \Omega_{1 \cdots (n-2)\alpha\beta} h_{\alpha(n-1)k} h_{\beta n k}]$$
(2.1)

*where* $\Omega_{\alpha\beta 3 \cdots n} = \Omega(e_\alpha, e_\beta, e_3, \cdots, e_n)$ *and*

$$(\frac{d}{dt} - \Delta)|A|^2 = -2|\nabla A|^2 + 2 \sum_{\alpha,\gamma,i,m} (\sum_k h_{\alpha i k} h_{\gamma m k} - h_{\alpha m k} h_{\gamma i k})^2 + 2 \sum_{i,j,m,k} (\sum_\alpha h_{\alpha i j} h_{\alpha m k})^2$$
(2.2)

*where $\Delta$ is the Laplace operator of the induced metric on $\Sigma_t$.*

*Proof.* We shall compute the first formula in terms of local coordinates. For a calculation of the first and second formula in moving frames, please see [7].

Let $F_t : \Sigma \to \mathbb{R}^{n+m}$ be a mean curvature flow given by

$$F_t = (F^1(x^1, \cdots, x^n, t), \cdots, F^N(x^1, \cdots, x^n, t))$$

where $(x_1, \cdots, x^n)$ are local coordinates on $\Sigma$. In the following calculation, we shall suppress the time index $t$.

The induced metric is $g_{ij} = \langle \frac{\partial F}{\partial x^i}, \frac{\partial F}{\partial x^j} \rangle = I_{AB} \frac{\partial F^A}{\partial x^i} \frac{\partial F^B}{\partial x^j}$.

Let $\Omega$ be a parallel $n$-form on $\mathbb{R}^{n+m}$. We are interested in the time-dependent function

$$*\Omega = \frac{\Omega(\frac{\partial F}{\partial x^1}, \cdots, \frac{\partial F}{\partial x^n})}{\sqrt{\det g_{ij}}}.$$

Now

$$\frac{d}{dt} * \Omega = \frac{1}{\sqrt{\det g_{ij}}} \frac{d}{dt}(\Omega(\frac{\partial F}{\partial x^1}, \cdots, \frac{\partial F}{\partial x^n})) + *\Omega \frac{d}{dt} \frac{1}{\sqrt{\det g_{ij}}}$$

and

$$\frac{d}{dt}(\Omega(\frac{\partial F}{\partial x^1}, \cdots, \frac{\partial F}{\partial x^n})) = \Omega(\frac{\partial H}{\partial x^1} \cdots \frac{\partial F}{\partial x^n}) + \cdots + \Omega(\frac{\partial F}{\partial x^1} \cdots \frac{\partial H}{\partial x^n}).$$

Recall along any mean curvature flow, we have

$$\frac{d}{dt}\sqrt{\det g_{ij}} = -|H|^2\sqrt{\det g_{ij}}$$

After this calculation we no longer need to vary the time variable and at a given time and a point $p$ we may assume $g_{ij} = \delta_{ij}$ and $\det g_{ij} = 1$. Therefore

$$\frac{d}{dt} *\Omega = \Omega(\frac{\partial H}{\partial x^1}\cdots \frac{\partial F}{\partial x^n}) + \cdots + \Omega(\frac{\partial F}{\partial x^1}\cdots \frac{\partial H}{\partial x^n}) + *\Omega|H|^2.$$

Now decompose $\frac{\partial H}{\partial x^i}$ into normal and tangent parts,

$$\frac{\partial H}{\partial x^i} = \nabla_i H + (\frac{\partial H}{\partial x^i})^T = \nabla_i H - \langle H, \frac{\partial^2 F}{\partial x^i \partial x^j}\rangle \frac{\partial F}{\partial x^j}$$

where $\nabla_i H$ is the covariant derivative of $H$ as a normal vector field.

Thus

$$\frac{d}{dt} *\Omega = \Omega(\nabla_1 H \cdots \frac{\partial F}{\partial x^n}) + \cdots \Omega(\frac{\partial F}{\partial x^1} \cdots \nabla_n H). \qquad (2.3)$$

At this point $p$, we may further assume we have a normal coordinate system so that $\frac{\partial}{\partial x^i}(g_{kl}) = 0$ at $p$ and thus $\sum_l g_{kl}\Gamma^l_{ij} = \langle \frac{\partial^2 F}{\partial x^k \partial x^i}, \frac{\partial F}{\partial x^j}\rangle = 0$. $\frac{\partial^2 F}{\partial x^k \partial x^i}$ is in the normal direction representing the second fundamental form and is denoted by $H_{ki}$.

Now the induced Laplacian at $p$ is $\Delta = g^{kl}\frac{\partial}{\partial x^k}\frac{\partial}{\partial x^l} = \frac{\partial^2}{\partial x^k \partial x^k}$. Therefore,

$$\Delta *\Omega = \frac{\partial^2}{\partial x^k \partial x^k}(\Omega(\frac{\partial F}{\partial x^1},\cdots,\frac{\partial F}{\partial x^n})) + *\Omega\frac{\partial^2}{\partial x^k \partial x^k}(\frac{1}{\sqrt{\det g_{ij}}}). \qquad (2.4)$$

The first term on the righthand side is

$$\frac{\partial}{\partial x^k}(\Omega(\frac{\partial}{\partial x^k}\frac{\partial F}{\partial x^1},\cdots,\frac{\partial F}{\partial x^n}) + \cdots + \Omega(\frac{\partial F}{\partial x^1},\cdots,\frac{\partial}{\partial x^k}\frac{\partial F}{\partial x^n}))$$
$$= \Omega(\frac{\partial^2}{\partial x^k \partial x^k}\frac{\partial F}{\partial x^1},\cdots\frac{\partial F}{\partial x^n}) + \cdots + \Omega(\frac{\partial F}{\partial x^1},\frac{\partial^2}{\partial x^k \partial x^k}\frac{\partial F}{\partial x^n})$$
$$+ 2[\Omega(H_{k1}, H_{k2}, \frac{\partial F}{\partial x^3},\cdots,\frac{\partial F}{\partial x^n}) + \cdots \Omega(\frac{\partial F}{\partial x^1},\cdots, H_{k\,n-1}, H_{k\,n})].$$

Decompose $\frac{\partial^3 F^A}{\partial x^k \partial x^k \partial x^i}$ into tangent and normal parts:

$$\frac{\partial^3 F^A}{\partial x^k \partial x^k \partial x^i} = (I_{AB}\frac{\partial^3 F^A}{\partial x^k \partial x^k \partial x^i}\frac{\partial F^B}{\partial x^m})\frac{\partial F^A}{\partial x^m} + (\frac{\partial^3 F^A}{\partial x^k \partial x^k \partial x^i})^\perp.$$

Now
$$\left(\frac{\partial^3 F^A}{\partial x^k \partial x^k \partial x^i}\right)^\perp$$
$$= \left[\frac{\partial}{\partial x^k}\left(\Gamma^p_{ki}\frac{\partial F}{\partial x^p} + H_{ki}\right)\right]^\perp$$
$$= \left(\frac{\partial}{\partial x^k}H_{ki}\right)^\perp$$
$$= \nabla_k H_{ki}$$
$$= \nabla_i H$$

where we use the Codazzi equation $\nabla_k H_{ij} = \nabla_i H_{kj}$.

Therefore, the first term on the right hand side of equation of (2.4) becomes

$$*\Omega I_{AB}\frac{\partial^3 F^A}{\partial x^k \partial x^k \partial x^i}\frac{\partial F^B}{\partial x^i}$$
$$+ \Omega(\nabla_1 H, \frac{\partial F}{\partial x^2}, \cdots, \frac{\partial F}{\partial x^n}) + \Omega(\frac{\partial F}{\partial x^1}, \cdots, \nabla_n H)$$
$$+ 2[\Omega(H_{k1}, H_{k2}, \frac{\partial F}{\partial x^3}, \cdots, \frac{\partial F}{\partial x^n}) + \cdots \Omega(\frac{\partial F}{\partial x^1}, \cdots, H_{k\,n-1}, H_{k\,n})].$$

The second term on the right hand side of equation (2.4) is

$$*\Omega\frac{\partial^2}{\partial x^k \partial x^k}\left(\frac{1}{\sqrt{\det g_{ij}}}\right) = *\Omega\frac{\partial}{\partial x^k}\left(-\frac{1}{2}(\det g_{ij})^{-3/2}\frac{\partial}{\partial x^k}(\det g_{ij})\right)$$
$$= *\Omega\frac{\partial}{\partial x^k}\left[-\frac{1}{2}(\det g_{ij})^{-3/2}\left(\frac{\partial}{\partial x^k}(g_{ij})g^{ji}\right)\det g_{ij}\right]$$
$$= -\frac{1}{2} * \Omega\left(\left(\frac{\partial^2}{\partial x^k \partial x^k}g_{ij}\right)g^{ji}\right).$$

In the last equality, we use $\sqrt{\det g_{ij}} = 1$ at $p$.

Now

$$\frac{\partial^2}{\partial x^k \partial x^k}g_{ij} = \frac{\partial^2}{\partial x^k \partial x^k}\left(I_{AB}\frac{\partial F^A}{\partial x^i}\frac{\partial F^B}{\partial x^j}\right)$$
$$= 2\frac{\partial}{\partial x^k}\left(I_{AB}\frac{\partial^2 F^A}{\partial x^k \partial x^i}\frac{\partial F^B}{\partial x^j}\right)$$
$$= 2\left(I_{AB}\frac{\partial^3 F^A}{\partial x^k \partial x^k \partial x^i} + I_{AB}\frac{\partial^2 F^A}{\partial x^k \partial x^i}\frac{\partial^2 F^B}{\partial x^k \partial x^j}\right)$$

Therefore the second term in equation (2.4) is

$$*\Omega\frac{\partial^2}{\partial x^k \partial x^k}(\frac{1}{\sqrt{\det g_{ij}}}) = *\Omega\sum_{i,k}[-I_{AB}\frac{\partial^3 F^A}{\partial x^k \partial x^k \partial x^i}\frac{\partial F^B}{\partial x^i} - I_{AB}\frac{\partial^2 F^A}{\partial x^k \partial x^i}\frac{\partial^2 F^B}{\partial x^k \partial x^i}]$$

We arrive at

$$\Delta * \Omega = - *\Omega |A|^2$$
$$+ \Omega(\nabla_1 H, \frac{\partial F}{\partial x^2}, \cdots, \frac{\partial F}{\partial x^n}) + \Omega(\frac{\partial F}{\partial x^1}, \cdots, \nabla_n H)$$
$$+ 2[\Omega(H_{k1}, H_{k2}, \frac{\partial F}{\partial x^3}, \cdots, \frac{\partial F}{\partial x^n}) + \cdots \Omega(\frac{\partial F}{\partial x^1}, \cdots, H_{k\,n-1}, H_{k\,n})].$$

Combine this with equation (2.3) we obtain the desired formula. □

The proof of Theorem A utilizes the evolution equation of $*\Omega$ when $\Omega$ is not necessarily parallel. The derivation is the same as the parallel case except the derivatives of $\Omega$ will be involved. The formula for a general ambient Riemannian manifold is derived in §3 of [9]. Let $y^1, \cdots, y^{n+m}$ be the fixed coordinates on $\mathbb{R}^{n+m}$ and $\Omega = \sum_{A_1 < \cdots < A_n} \Omega_{A_1, \cdots A_n} dy^1 \wedge \cdots \wedge dy^{A_n}$ be a general $n$ form in $\mathbb{R}^{n+m}$.

**Lemma 2.2** *Let $\Omega$ be a general $n$-form on $\mathbb{R}^{n+m}$. At a point $p$ of $\Sigma_t$, we choose our coordinates so that $\{\frac{\partial F}{\partial x^i}\}$ is orthonormal and $\frac{\partial^2 F}{\partial x^i \partial x^j}$ is in the normal direction. With respect to the orthonormal basis $\{e_i = \frac{\partial F}{\partial x^i}\}$ we have*

$$(\frac{d}{dt} - \Delta) * \Omega = *\Omega|A|^2 - 2\sum_{\alpha,\beta,k}[\Omega_{\alpha\beta 3\cdots n}h_{\alpha 1 k}h_{\beta 2 k} + \cdots + \Omega_{1\cdots(n-2)\alpha\beta}h_{\alpha(n-1)k}h_{\beta n k}]$$
$$- *(tr_{\Sigma_t} D^2\Omega) - 2\sum_{\alpha,k}[(D_{e_k}\Omega)_{\alpha 2\cdots n}h_{\alpha 1 k} + \cdots + (D_{e_k}\Omega)_{1\cdots n-1,\alpha}h_{\alpha n k}]$$
(2.5)

where

$$tr_{\Sigma_t} D^2\Omega = \sum_{A_1 < \cdots < A_n, k} (D^2\Omega_{A_1, \cdots, A_n})(e_k, e_k)dy^1 \wedge \cdots \wedge dy^{A_n}$$

and

$$D_{e_k}\Omega = \sum_{A_1 < \cdots < A_n} \langle D\Omega_{A_1, \cdots, A_n}, e_k \rangle dy^1 \wedge \cdots \wedge dy^{A_n}.$$

$D^2\Omega_{A_1,\cdots,A_n}$ and $D\Omega_{A_1,\cdots,A_n}$ are just the ordinary Hessian and gradient of the function $\Omega_{A_1,\cdots,A_n}$ on $\mathbb{R}^{n+m}$.

*Proof.* This follows from the general formula in §3 of [9] by noting that, in the notation of [9],
$$(\nabla^M)^2(e_k, e_k)\Omega = \nabla^M_{e_k}\nabla^M_{e_k}\Omega - \nabla^M_{\nabla^M_{e_k}e_k}\Omega$$
and $\nabla^M_{e_k}e_k = H$ by our choice of $(\nabla^M_{e_k}e_k)^T = 0$ at the point of calculation. □

Let $u$ denote the distance function to the reference $\mathbb{R}^n$, then $u^2 = \sum_\alpha \langle F, a_\alpha\rangle^2$. The following two equations will be used in the localization of estimates.

**Lemma 2.3**
$$(\frac{d}{dt} - \Delta)|F|^2 = -2n \tag{2.6}$$

$$(\frac{d}{dt} - \Delta)u^2 = -2\sum_{i,\alpha}\langle e_i, a_\alpha\rangle^2 \tag{2.7}$$

*Proof.* The equation for $|F|^2$ is derived as the following:
$$\frac{d}{dt}|F|^2 = 2\langle\frac{dF}{dt}, F\rangle = 2\langle\Delta F, F\rangle,$$

$$\Delta|F|^2 = 2\langle\Delta F, F\rangle + 2|\nabla F|^2,$$

and
$$|\nabla F|^2 = \sum_i |\nabla_{e_i}F|^2 = \sum_i |e_i|^2 = n.$$

The equation for $u^2$ is derived similarly.
$$\frac{d}{dt}u^2 = 2\sum_\alpha\langle F, a_\alpha\rangle\langle\frac{dF}{dt}, a_\alpha\rangle = 2\sum_\alpha\langle F, a_\alpha\rangle\langle\Delta F, a_\alpha\rangle$$

where we used $H = \Delta_\Sigma F$. Now
$$\Delta u^2 = 2\sum_\alpha\langle\nabla F, a_\alpha\rangle^2 + 2\langle F, a_\alpha\rangle\langle\Delta F, a_\alpha\rangle.$$

The derivation is completed by noting $\sum_\alpha\langle\nabla F, a_\alpha\rangle^2 = \sum_\alpha\sum_i\langle\nabla_{e_i}F, a_\alpha\rangle$.
□

Given any submanifold $\Sigma$ of $\mathbb{R}^{n+m}$, denote by $N_\delta(\Sigma) = \{z|d(z,\Sigma) < \delta\}$ the $\delta$ tubular neighborhood of $\Sigma$ in $\mathbb{R}^{n+m}$. The following proposition controls the Hausdorff distance between $\Sigma$ and $\Sigma_t$ along the mean curvature flow.

**Proposition 2.1** *Let $\Sigma$ be a smooth compact submanifold in $\mathbb{R}^{n+m}$, then*

$$N_{\sqrt{\delta^2-2nt}}(\Sigma_t) \subset N_\delta(\Sigma).$$

*Proof.* We claim for $y_0 \in \mathbb{R}^{n+m}$ if $B_\delta(y_0) \cap \Sigma$ is empty, so is $B_{\sqrt{\delta^2-2nt}}(y_0) \cap \Sigma_t$. We may assume $y_0$ is the origin, let

$$\phi = |F(x,t)|^2 - \delta^2 + 2nt.$$

By equation (2.6), $(\frac{d}{dt} - \Delta)\phi = 0$. $\phi > 0$ at $t = 0$, so by the maximum principle, we have $\phi > 0$ afterwards. □

Set $\delta = \sqrt{2nt}$, we obtain

**Corollary 2.1** $\Sigma_t$ *remains in the $\sqrt{2nt}$ neighborhood of $\Sigma$.*

## 3 Local gradient estimates

If $\Sigma_t$ is locally given as the graph of a vector-valued function $f : U \subset \mathbb{R}^n \to \mathbb{R}^m$. $*\Omega$ is in fact the Jacobian of the projection map from $\Sigma_t$ to $\mathbb{R}^n$ and in terms of $f$

$$*\Omega = \frac{1}{\sqrt{\det(I + (df)^T df)}}$$

where $(df)^T$ is the adjoint of $df$. Any lower bound of $*\Omega$ gives an upper bound for $|df|$.

**Lemma 3.1** *If $*\Omega > \frac{1}{\sqrt{2}}$, then*

$$(\frac{d}{dt} - \Delta) *\Omega \geq (2 - \frac{1}{(*\Omega)^2})|A|^2.$$

*Proof.* As in [8], we can rewrite equation (2.1) in terms of the singular values $\lambda_i, i = 1 \cdots, \min\{n,m\}$ of $df$, for any local defining function of $\Sigma_t$.

$$(\frac{d}{dt} - \Delta) *\Omega = *\Omega\{\sum_{\alpha,l,k} h^2_{\alpha lk} - 2\sum_{k,i<j} \lambda_i\lambda_j h_{n+i,ik}h_{n+j,jk} + 2\sum_{k,i<j} \lambda_i\lambda_j h_{n+j,ik}h_{n+i,jk}]\}$$

(3.1)

where the index $i, j$ runs from 1 to $\min\{n, m\}$.

Because $*\Omega - \frac{1}{\sqrt{2}} > 0$ and $*\Omega = \frac{1}{\sqrt{\prod_{i=1}^n(1+\lambda_i^2)}}$, we have $\prod_{i=1}^n(1+\lambda_i^2) < 2$.
Denote $\delta = 2 - \prod_{i=1}^n(1+\lambda_i^2) > 0$. It is clear that

$$\prod_i(1+\lambda_i^2) \geq 1 + \sum_i \lambda_i^2$$

and we derive

$$\sum_i \lambda_i^2 \leq 1 - \delta, \text{ and } |\lambda_i \lambda_j| \leq 1 - \delta$$

As in [8] by completing square we obtain

$$(\frac{d}{dt} - \Delta) * \Omega \geq \delta |A|^2.$$

□

Following the notation in Theorem 2.1 of [4]. Let $y_0$ be an arbitrary point in $\mathbb{R}^{n+m}$. $\phi(y,t) = R^2 - |y - y_0|^2 - 2nt$ and $\phi_+$ denotes the positive part of $\phi$.

**Lemma 3.2** *Let $f$ be a function defined on $\Sigma_t$ with $0 < f < K$ such that $(\frac{d}{dt} - \Delta)f \geq 0$. Let $v = \frac{1}{f}$ then $v$ satisfies*

$$v(F,t)\phi_+(F,t) \leq \sup_{\Sigma_0} v\phi_+$$

*Proof.* The proof is adapted from Theorem 2.1 in [4]. We may assume $y_0$ is the origin. Set

$$\eta(x,t) = (|F(x,t)|^2 - R^2 + 2nt)^2.$$

Then by equation (2.6),

$$(\frac{d}{dt} - \Delta)\eta = -2|\nabla |F|^2|^2.$$

Now

$$(\frac{d}{dt} - \Delta)v^2 = -2f^{-3}(\frac{d}{dt} - \Delta)f - 6f^{-4}|\nabla f|^2 \leq -6|\nabla v|^2.$$

Combine these equations we obtain

$$(\frac{d}{dt}-\Delta)(v^2\eta) = (\frac{d}{dt}-\Delta)v^2\eta+(\frac{d}{dt}-\Delta)\eta v^2-2\nabla v^2\cdot\nabla\eta \leq -6|\nabla v|^2\eta-2|\nabla|F|^2|^2v^2-2\nabla v^2\cdot\nabla\eta.$$

As in [4], observe that

$$-2\nabla v^2 \cdot \nabla\eta = -6v\nabla v\nabla\eta + \eta^{-1}\nabla\eta\nabla(v^2\eta) - 4|\nabla|F|^2|^2v^2.$$

Therefore

$$(\frac{d}{dt} - \Delta)(v^2\eta) \leq -6|\nabla v|^2\eta - 6|\nabla|F|^2|^2v^2 - 6v\nabla v \cdot \nabla\eta + \eta^{-1}\nabla\eta \cdot \nabla(v^2\eta).$$

Apply Young's inequality to the term $6v\nabla v \cdot \nabla\eta$.

$$6v\nabla v \cdot \nabla\eta \leq 6|\nabla v|^2\eta + \frac{3}{2}v^2\eta^{-1}(\eta')^2|\nabla|F|^2|^2 = 6|\nabla v|^2\eta + 6|\nabla|F|^2|^2v^2$$

The estimate follows by replacing $\eta$ by $\phi_+$ and applying the maximum principle. □

## 4  Local curvature estimates

The following local comparison lemma generalizes Theorem 3.1 in [4]. It applies to other heat equations and is interesting in its own right. Let $r$ be a non-negative function on $\Sigma_t$ satisfying

$$|(\frac{d}{dt} - \Delta)r| \leq c_4$$

and

$$|\nabla r|^2 \leq c_5 r.$$

In later application $r(x,t)$ will be $|F(x,t)-y_0|^2+2nt$ or $|F(x,t)|^2-u^2(F(x,t))$. Both functions satisfy the above assumptions by equations (2.6) and (2.7).

**Lemma 4.1** *If $h$ and $f \leq c_3$ are positive functions on $\Sigma_t$ that satisfy*

$$(\frac{d}{dt} - \Delta)h \leq c_1 h^3$$

and
$$(\frac{d}{dt} - \Delta)f \geq c_2 f h^2.$$

Let $R > 0$ be a constant such that $\Sigma_t \cap \{r(x,t) \leq R^2\}$ is compact and $0 < \theta < 1$. If $0 < c_1 \leq c_2$ then the following estimate holds

$$\sup_{r(x,t) \leq \theta R^2} h^2 \leq c_0 (1-\theta)^{-2}(t^{-1} + R^{-2}) \sup_{s \in [0,t]} \sup_{r(x,s) \leq R^2} \frac{1}{f^4}$$

where $c_0$ is a constant depending only on $c_1$, $c_3$, $c_4$ and $c_5$.

*Proof.* First of all, we calculate the evolution equation satisfied by $h^2$:

$$(\frac{d}{dt} - \Delta)h^2 = 2h(\frac{d}{dt} - \Delta)h - 2|\nabla h|^2 \leq 2c_1 h^4 - 2|\nabla h|^2.$$

Denote $v = \frac{1}{f}$, then

$$(\frac{d}{dt} - \Delta)v^2 = -2f^{-3}(\frac{d}{dt} - \Delta)f - 6f^{-4}|\nabla f|^2 \leq -2c_2 v^2 h^2 - 6|\nabla v|^2.$$

As in [4], we consider $\phi(v^2)$ for a positive and increasing function $\phi$ to be determined later. Thus

$$(\frac{d}{dt} - \Delta)\phi = \phi'(\frac{d}{dt} - \Delta)v^2 - \phi''|\nabla v^2|^2 \leq -2c_2 \phi' v^2 h^2 - 6\phi'|\nabla v|^2 - 4\phi'' v^2 |\nabla v|^2.$$

We obtain

$$(\frac{d}{dt} - \Delta)(h^2 \phi) = (\frac{d}{dt} - \Delta)h^2 \phi + (\frac{d}{dt} - \Delta)\phi h^2 - 2\nabla h^2 \cdot \nabla \phi$$
$$\leq 2(c_1 \phi - c_2 \phi' v^2)h^4 - 2|\nabla h|^2 \phi - (6\phi' + 4\phi'' v^2)h^2 |\nabla v|^2 - 2\nabla h^2 \cdot \nabla \phi.$$

Now break the term $-2\nabla h^2 \cdot \nabla \phi$ into two equal parts. One of them is

$$-\nabla h^2 \cdot \nabla \phi = -\phi^{-1} \nabla \phi \nabla(h^2 \phi) + \phi^{-1}|\nabla \phi|^2 h^2,$$

while the other can be estimated by

$$-\nabla h^2 \cdot \nabla \phi = -2h \nabla h \cdot \nabla \phi \leq 2|\nabla h|^2 \phi + \frac{1}{2}\phi^{-1}|\nabla \phi|^2 h^2.$$

Since $\phi^{-1}|\nabla\phi|^2 h^2 = 4\phi^{-1}(\phi')^2 v^2 |\nabla v|^2 h^2$, we arrive at

$$(\frac{d}{dt}-\Delta)(h^2\phi) \leq 2(c_1\phi - c_2\phi'v^2)h^4 - \phi^{-1}\nabla\phi \cdot \nabla(h^2\phi) - (6\phi'(1-\phi^{-1}\phi'v^2) + 4\phi''v^2)h^2|\nabla v|^2.$$

Now using an idea of Caffarelli, Nirenberg and Spruck in [2], we set $\phi(v^2) = \frac{v^2}{1-kv^2}$ where $k = \frac{1}{2}\inf_{s\in[0,t]}\inf_{\{r(x,s)\leq R^2\}} v^{-2}$. Then

$$c_1\phi - c_2\phi'v^2 = \frac{(c_1-c_2)v^2 - c_1kv^4}{(1-kv^2)^2} \leq -c_1k\phi^2,$$

$$6\phi'(1-\phi^{-1}\phi'v^2) + 4\phi''v^2 = \frac{2k}{(1-kv^2)^2}\phi,$$

and

$$\phi^{-1}\nabla\phi = 2\phi v^{-3}\nabla v.$$

Thus $g = h^2\phi$ satisfies

$$(\frac{d}{dt}-\Delta)g \leq -2c_1kg^2 - \frac{2k}{(1-kv^2)^2}|\nabla v|^2 g - 2\phi v^{-3}\nabla v \cdot \nabla g.$$

The rest of the proof is the same as that of Theorem 3.1 in [4]. The term $-\frac{2k}{(1-kv^2)^2}|\nabla v|^2 g$ helps to cancel similar order terms in later calculations. $\square$

**Corollary 4.1** *If $P$ is a positive bounded function satisfying*

$$(\frac{d}{dt}-\Delta)P \geq 5P|A|^2$$

*on $\Sigma_s \cap \{|y-y_0| \leq R^2 - 2ns\}$, $s \in [0,t]$, then for any $0 < \theta < 1$,*

$$\sup_{\Sigma_t \cap \{|y-y_0|^2 \leq \theta R^2 - 2nt\}} |A|^2 \leq c_0(1-\theta)^{-2}(\frac{1}{R^2} + \frac{1}{t}) \sup_{s\in[0,t]} \sup_{\Sigma_s \cap \{|y-y_0|^2 \leq R^2 - 2ns\}} \frac{1}{P^4}.$$

*Proof.* We first estimate the higher order terms in equation (2.2):

$$\sum_{\alpha,\gamma,i,m}(\sum_k h_{\alpha ik}h_{\gamma mk} - h_{\alpha mk}h_{\gamma ik})^2$$

$$\leq \sum_{\alpha,\gamma,i,m}[2(\sum_k h_{\alpha ik}h_{\gamma mk})^2 + 2(\sum_k h_{\alpha mk}h_{\gamma ik})^2]$$

$$\leq 2\sum_{\alpha,\gamma,i,m}(\sum_k h_{\alpha ik}^2)(\sum_k h_{\gamma mk}^2) + 2\sum_{\alpha,\gamma,i,m}(\sum_k h_{\alpha mk}^2)(\sum_k h_{\gamma ik}^2)$$

$$\leq 4|A|^4$$

and

$$\sum_{i,j,m,k}(\sum_\alpha h_{\alpha ij}h_{\alpha mk})^2$$
$$\sum_{i,j,m,k}(\sum_\alpha h_{\alpha ij}^2)(\sum_\alpha h_{\alpha mk}^2)$$
$$\leq |A|^4.$$

Therefore we have

$$(\frac{d}{dt}-\Delta)|A|^2 \leq 10|A|^4 - 2|\nabla A|^2.$$

Since

$$(\frac{d}{dt}-\Delta)|A|^2 = 2|A|(\frac{d}{dt}-\Delta)|A| - 2|\nabla|A||^2$$

Schwartz' inequality gives $|\nabla|A||^2 \leq |\nabla A|^2$. Therefore

$$(\frac{d}{dt}-\Delta)|A| \leq 5|A|^3$$

The corollary follows by choosing $r = |F(x,t) - y_0|^2 + 2nt$ in Lemma 4.1. □

## 5  Proof of Theorem A

First we define the local Lipschitz condition:

**Definition 5.1** *Given any positive $K < 1$, a compact n-dimensional submanifold $\Sigma$ of $\mathbb{R}^{n+m}$ is said to satisfy the $K$ local Lipschitz condition if there exists a $r_0 >$ such that $\Sigma \cap B_q(r_0)$ for each $q \in \Sigma$ can be written as the graph of a vector valued Lipschitz function $f_q$ over an n-dimensional affine space $L_q$ through $q$ with $\frac{1}{\sqrt{\det(I+(df_q)^T df_q)}} > K$.*

We are ready to prove Theorem A.

*Proof of Theorem A.* First we assume $\Sigma$ is a smooth compact submanifold that satisfies the $K$ local Lipschitz condition for $K$ to be determined in the proof. Denote the volume form of $L_q$ by $\Omega_{L_q}$, the assumption implies $*\Omega_{L_q} > K$ on $\Sigma \cap B_q(r_0)$. The idea is to construct an $n$-form $\Omega$ by "averaging" $\Omega_{L_q}$ so that $*\Omega$ gives the desired positive function in Corollary 4.1.

Choose a subset $\{q_\nu\}_{\nu \in \Lambda} \subset \Sigma$ such that $\{B_{q_\nu}(\frac{r_0}{5})\}_{\nu \in \Lambda}$ are maximally pairwise disjoint and $\cup_{q \in \Sigma} B_q(\frac{r_0}{5}) \subset \cup_{\nu \in \Lambda} B_{q_\nu}(r_0)$. We can arrange that for each $q \in \Sigma$, there exists a $q_\nu$ with $B_q(\frac{r_0}{5}) \subset B_{q_\nu}(\frac{4r_0}{5})$.

It is clear that
$$N_{\frac{r_0}{5}}(\Sigma) \subset \cup_{q \in \Sigma} B_q(\frac{r_0}{5}).$$

Set $t_0 = \frac{r_0^2}{50n}$, by Corollary 2.1 we have

$$\Sigma_t \subset N_{\frac{r_0}{5}}(\Sigma) \text{ for } t \in [0, t_0]. \tag{5.1}$$

Let $\Omega_\nu$ be the volume form of the affine space $L_{q_\nu}$ and $\phi_\nu$ be a cut-off function such that $0 < \phi_\nu < 1$ on $B_{q_\nu}(r_0)$, $\phi_\nu = 1$ on $B_{q_\nu}(\frac{4r_0}{5})$ and $\phi_\nu = 0$ outside $B_{q_\nu}(r_0)$. Therefore $|D\phi_\nu| \leq \frac{c_1}{r_0}$ and $|D^2\phi_\nu| \leq \frac{c_2}{r_0^2}$. Take $p_\mu = \frac{\phi_\mu}{\sum_{\nu \in \Lambda} \phi_\nu}$ to be the partition of unity of $\cup_{\nu \in \Lambda} B_{q_\nu}(r_0)$. We claim

**Lemma 5.1** *For any $\mu \in \Lambda$ and $y \in N_{\frac{r_0}{5}}(\Sigma)$, we have $|Dp_\mu|(y) \leq \frac{c_3}{r_0}$ and $|D^2 p_\mu|(y) \leq \frac{c_4}{r_0^2}$ where $c_3$ and $c_4$ depend on $n$ and $m$.*

*Proof.* Given any $y \in N_{\frac{r_0}{5}}(\Sigma)$, there is a $q \in \Sigma$ and $q_\nu$ so that $y \in B_q(\frac{r_0}{5}) \subset B_{q_\nu}(\frac{4r_0}{5})$. Therefore $\phi_\nu(y) = 1$ for this $\nu$ and thus $\sum_{\nu \in \Lambda} \phi_\nu(y) > 1$ for any $y \in N_{\frac{r_0}{5}}(\Sigma)$. On the other hand, the number of $\phi_\nu$ such that $\phi_\nu(y) > 0$ for any $y$ is bounded by a constant depending on $n + m$.

The lemma follows by estimating
$$Dp_\mu = \frac{D\phi_\mu}{\sum_\nu \phi_\nu} - \frac{\phi_\mu D(\sum_\nu \phi_\nu)}{(\sum_\nu \phi_\nu)^2}$$

and

$$D^2 p_\mu = \frac{D^2 \phi_\mu}{\sum_\nu \phi_\nu} - 2\frac{(D\phi_\mu)D(\sum_\nu \phi_\nu)}{(\sum_\nu \phi_\nu)^2} - \frac{\phi_\mu D^2(\sum_\nu \phi_\nu)}{(\sum_\nu \phi_\nu)^2} + 2\frac{\phi_\mu |D(\sum_\nu \phi_\nu)|^2}{(\sum_\nu \phi_\nu)^3}.$$

$\square$

Now define
$$\Omega = \sum_{\nu \in \Lambda} p_\nu \Omega_\nu, \tag{5.2}$$

$\Omega$ is an $n$-form on $\mathbb{R}^{n+m}$ supported in $\cup_{\nu \in \Lambda} B_{q_\nu}(r_0)$.

We notice that $D\Omega = \sum_\nu (Dp_\nu)\Omega_\nu$ and $D^2\Omega = \sum_\nu (D^2 p_\nu)\Omega_\nu$, therefore by Lemma 2.2, we derive

$$(\frac{d}{dt} - \Delta) * \Omega \geq *\Omega |A|^2 - 2 \sum_{\alpha,\beta,k} [\Omega_{\alpha\beta 3 \cdots n} h_{\alpha 1 k} h_{\beta 2 k} + \cdots + \Omega_{1 \cdots (n-2)\alpha\beta} h_{\alpha(n-1)k} h_{\beta n k}]$$
$$- c_5 - c_6 |A|$$
(5.3)

where $c_5$ and $c_6$ are constants depending on $n, m$, and $r_0$.

Initially $*\Omega_\nu > K$ whenever $p_\nu > 0$ and $\sum_{p_\nu > 0} p_\nu = 1$, so we have $*\Omega = \sum_\nu p_\nu * \Omega_\nu > K$ on $\Sigma$. Because each $\Omega_\nu$, as a vector in $\wedge^n (\mathbb{R}^{n+m})^*$, has $L^2$ norm $|\Omega_\nu|^2 = 1$, we have $|\Omega|^2 \leq 1$ as $\Omega$ is a convex combination of $L^2$ unit vectors.

Let $K_0 < K$ be a positive constant to be determined and by assumption $*\Omega > K > K_0$ on $\Sigma$. For any positive constant $\epsilon < K - K_0$ (for example $\epsilon = \frac{K-K_0}{2}$), set $c_7 = c_5 + \frac{1}{4\epsilon} c_6^2$. We claim

$$*\Omega + c_7 t > K \tag{5.4}$$

on $\Sigma_t$ for $t \in [0, t_1]$ where

$$t_1 = \min\{\frac{K - K_0 - \epsilon}{c_7}, t_0\}.$$

This is proved by the maximum principle. Suppose inequality (5.4) is violated at some $t_2 < t_1$ for the first time. Then $*\Omega + c_7 t \geq K$ on $[0, t_2]$ and thus $*\Omega \geq K_0 + \epsilon$ on $[0, t_2]$. This implies $\Omega_{\alpha\beta 3 \cdots n} \leq \sqrt{1 - K_0^2}$ since $e_1 \wedge \cdots \wedge e_n$ is orthogonal to $e_\alpha \wedge e_\beta \wedge e_3 \wedge \cdots \wedge e_n$ and $|\Omega|^2 \leq 1$. This estimate holds for other similar terms inside the bracket on the right hand side of inequality (5.3) . It follows that

$$(\frac{d}{dt} - \Delta) * \Omega \geq (*\Omega - c_8 \sqrt{1 - K_0^2})|A|^2 - c_5 - c_6|A| \text{ on } [0, t_2].$$

where $c_8$ is a combinatorial constant depending on $n$ and $m$.

Now we choose $K_0$ close to 1 so that $K_0 - c_8\sqrt{1 - K_0^2} = c_9 > 0$. Again $c_9$ depends on $n$ and $m$. Using the inequality $c_6|A| \leq \frac{1}{4\epsilon} c_6^2 + \epsilon |A|^2$, we obtain

$$(\frac{d}{dt} - \Delta) * \Omega \geq (*\Omega - c_8\sqrt{1 - K_0^2} - \epsilon)|A|^2 - c_5 - \frac{1}{4\epsilon} c_6^2$$

or
$$(\frac{d}{dt} - \Delta)(*\Omega + c_7 t) \geq (*\Omega - c_8\sqrt{1 - K_0^2} - \epsilon)|A|^2. \tag{5.5}$$

Because $*\Omega > K_0 + \epsilon$, it is easy to see $*\Omega - c_8\sqrt{1 - K_0^2} - \epsilon > K_0 - c_8\sqrt{1 - K_0^2} = c_9$. By the maximum principle $*\Omega + c_7 t$ is increasing in $[0, t_1]$.

Next we claim we can choose $K$ close to 1 to produce the desired positive bounded function in Corollary 4.1. Indeed, write equation (5.5) as

$$(\frac{d}{dt} - \Delta)[*\Omega + c_7 t - K] \geq \frac{*\Omega - c_8\sqrt{1 - K_0^2} - \epsilon}{*\Omega + c_7 t - K}[*\Omega + c_7 t - K]|A|^2$$

and recall this holds for $t \in [0, t_1]$. We claim if $K$ is close to 1, there exists a $T \leq t_1$ such that

$$\min_{0 \leq t \leq T} \min_{1 \geq x > K - c_7 t} \frac{x - c_8\sqrt{1 - K_0^2} - \epsilon}{x + c_7 t - K} > 5.$$

First of all, under $1 \geq x > K - c_7 t$, we have

$$\frac{x - c_8\sqrt{1 - K_0^2} - \epsilon}{x + c_7 t - K} > \frac{K - c_7 t - c_8\sqrt{1 - K_0^2} - \epsilon}{1 + c_7 t - K}$$

When $t = 0$ this is

$$\frac{K - c_8\sqrt{1 - K_0^2} - \epsilon}{1 - K}.$$

Since $K > K_0 + \epsilon$, the last expression is greater than

$$\frac{K_0 - c_8\sqrt{1 - K_0^2}}{1 - K} = \frac{c_9}{1 - K}.$$

Now we can choose $K$ close to 1 so that this is greater than 5. The dependence of the constants are as the following: $c_8$ is a combinatorial constant depending on $n$ and $m$, we determine $K_0$ and $c_9$ from

$$K_0 - c_8\sqrt{1 - K_0^2} = c_9$$

and then $K$ is determined by

$$\frac{c_9}{1 - K} > 5.$$

Therefore $K$ depends only on $n$ and $m$. The time $T$ is determined by

$$\frac{(K - c_7 t - c_8 \sqrt{1 - K_0^2} - \epsilon)}{1 + c_7 t - K} > 0$$

for $t \in [0, T]$ and $T \leq t_1 = \min\{\frac{K - K_0 - \epsilon}{c_7}, t_0\}$.

We then approximate a Lipschitz submanifold that satisfies the $K$ local Lipschitz condition by smooth submanifolds with the same bound on $*\Omega$. By Corollary 4.1, we have the uniform bound of $|A|^2$ and all higher derivatives bound for $|A|^2$ can be obtained as in [4]. Therefore the approximating mean curvature flows converge to a mean curvature flow that is smooth for $t \in (0, T]$ and the theorem is proved.

$\square$

We remark the theorem is also true when the ambient space is replaced by a complete Riemannian manifolds with bounded geometry. The ambient curvature only results in lower order terms in the evolution equation and the proof works if we shrink the time interval a little bit.